\documentclass
{article}
\usepackage{amsfonts, amsmath, amstext, amssymb}

\usepackage[dvips]{graphicx}
\usepackage{latexsym}

\newtheorem{thm}{Theorem}[section]
\newtheorem{prop}[thm]{Proposition}

\newtheorem{lemma}[thm]{Lemma}

\newtheorem{cor}[thm]{Corollary}

\newtheorem{definitiontemp}[thm]{Definition}
\newenvironment{defn}{\begin{definitiontemp}
\normalfont}{\end{definitiontemp}}

\def\bec{\begin{cor}}
\def\enc{\end{cor}}

\def\bet{\begin{thm}}
\def\ent{\end{thm}}

\def\becor{\begin{cor}}
\def\encor{\end{cor}}

\def\bel{\begin{lem}}
\def\enl{\end{lem}}

\def\bedef{\begin{defn}}
\def\endef{\end{defn}}

\def\bep{\begin{prop}}
\def\enp{\end{prop}}

\newenvironment{pf}{\begin{trivlist}\item[\hskip\labelsep
{\it Proof.}]}{\end{trivlist}}

\newcommand{\RCF}{\ensuremath{\textbf{RCF}}}

\newcommand{\Fin}{\textbf{Fin}}

\newcommand{\fin}{\texttt{fin}}
\newcommand{\Inf}{\textbf{Inf}}
\newcommand{\Ded}{\textbf{Ded}}
\newcommand{\Tot}{\textbf{Tot}}

\newcommand{\set}[2]{\ensuremath{ \{ #1 : #2 \} }}

\newcommand{\spec}[1]{\ensuremath{ \text{Spec} (#1 )}}

\renewcommand{\deg}[1]{\ensuremath{\text{deg}(#1)}}

\newcommand{\N}{\mathbb{N}}
\newcommand{\Z}{\mathbb{Z}}

\newcommand{\Q}{\mathbb{Q}}
\newcommand{\R}{\mathbb{R}}

\newcommand{\A}{\mathcal{A}}
\newcommand{\B}{\mathcal{B}}

\newcommand{\calS}{\mathcal{S}}

\newcommand{\Cbar}{\overline{C}}

\renewcommand{\L}{\mathcal{L}}

\renewcommand{\P}{\mathcal{P}}

\newcommand{\avec}{\vec{a}}

\newcommand{\Rc}{\R_{\bfc}}

\newcommand{\Rz}{\R_{\bfz}}

\newcommand{\la}{\langle}
\newcommand{\ra}{\rangle}

\newcommand{\comment}[1]{}

\def\diverges{\!\uparrow}
\def\converges{\!\downarrow}
\newcommand{\at}{\char'100}

\def\Ss{{\mathfrak S}}

\newcommand{\qed}{\hbox to 0pt{}\nobreak\hfill\rule{2mm}{2mm}}

\newcommand{\bfa}{\boldsymbol{a}}

\newcommand{\bfc}{\boldsymbol{c}}
\newcommand{\bfd}{\boldsymbol{d}}
\def\bfz{\boldsymbol{0}}
\def\s01{\ensuremath{\Sigma^0_1}}
\def\d02{\ensuremath{\Delta^0_2}}
\def\phi{\varphi}

\def\res{\!\!\upharpoonright\!}

\def\ep{\varepsilon}

\begin{document}

\title{Degree Spectra of Real Closed Fields}
\author{Russell Miller\thanks{This is a pre-print of an article to be published in the
\emph{Archive for Mathematical Logic}. The final
authenticated version will be available online at: https://doi.org/10.1007/s00153-018-0638-z.
The first author was
partially supported by Grant \# DMS -- 1362206 from
the National Science Foundation and by several grants from
the City University of New York PSC-CUNY Research Award Program.
The authors wish to acknowledge useful conversations
with Julia Knight and Reed Solomon.}
\ \& Victor Ocasio Gonz\'alez
}
\maketitle

\begin{abstract}
Several researchers have recently established that for every Turing degree
$\bfc$, the real closed field of all $\bfc$-computable
real numbers has spectrum $\set{\bfd}{\bfd'\geq\bfc''}$.
We investigate the spectra of real closed fields further,
focusing first on subfields of the field $\Rz$ of computable real numbers,
then on archimedean real closed fields more generally,
and finally on non-archimedean real closed fields.
For each noncomputable, computably enumerable set
$C$, we produce a real closed $C$-computable subfield of $\Rz$
with no computable copy.  Then we build an archimedean real closed field
with no computable copy but with a computable
enumeration of the Dedekind cuts it realizes, and a computably presentable
nonarchimedean real closed field whose residue field
has no computable presentation.
\end{abstract}

\section{Introduction}
\label{sec:intro}

The Turing degree spectrum of a countable first-order structure $\A$
provides a natural measure of the complexity of the isomorphism type
of that structure.  The spectrum of $\A$, by definition, is the set of
those Turing degrees $\bfd$ such that for some \emph{copy} of $\A$
(that is, for some $\B\cong\A$ with domain $\omega$), the atomic diagram
of $\B$ has Turing degree $\bfd$.  In \cite{K86}, Knight proved that the
spectrum is closed upwards under Turing reducibility in all but the most
trivial cases (in which it is a singleton), and so spectra fall under the more
general heading of \emph{mass problems}, a source of broad interest in
computability theory.

Nontrivial structures with computable copies (also known as computable
\emph{presentations}) have all Turing degrees in their spectrum, by Knight's result,
and we view these structures as being as simple as possible to present.
More complex structures have no computable copy, and the spectrum
tells us how much information suffices if one wishes to present a copy
of such a structure.  In certain cases, the spectrum may be the upper cone
above a specific degree $\bfd$, in which case we view $\bfd$ itself
as a very precise measure of the difficulty of presenting the structure,
as first described by Richter in \cite{R81}.  More recent work has yielded
many structures in which information (specifically, a given set $C$%
) is coded not directly into the
atomic diagram of the structure, but rather into the $n$-quantifier diagram,
for various $n$, and so the spectrum may contain those degrees whose
$n$-th jump computes $C$.  Such structures are discussed in
\cite{FHKKM12}, among other places, and we will meet some of them below.

Reversing these arguments, many researchers have also compared different
classes of structures by investigating the spectra which can be realized
by those classes.  The class of symmetric irreflexive graphs has been known
for quite some time to be complete in this sense:  every spectrum of a non-trivial
first-order structure (in a computable language) is also the spectrum
of such a graph.  Details appear in \cite{HKSS02}, which also established
this completeness for classes such as groups, rings, partial orders, and lattices.
Recently, in \cite{MPSS16}, fields of arbitrary characteristic have also been shown
to be complete in this sense.  On the other hand, classes such as
linear orders, Boolean algebras, and trees (as partial orders) were shown
in the same work of Richter not to be complete this way, and these classes have been parsed
further since then.  In \cite{DJ94}, for example, Boolean algebras were shown incapable
of realizing certain spectra which had been seen in \cite{JS91} to be the spectra
of linear orders.

To our knowledge, it remains unknown whether real closed fields have this completeness
property or not.  A real closed field is simply a model of the theory of the structure
$(\R,0,1,+,\cdot)$, the field of real numbers.  It makes no difference whether
we include the relation $<$ in the signature, since this relation is $\Delta^0_1$-definable
in every real closed field anyway:  $a < b$ if and only if $(b-a)$ has a nonzero square root
in the field, which holds if and only if $(a-b)$ has no square root.
Real closed fields have far more complex possible spectra than algebraically closed
fields, all of which have computable presentations.  Recently, independent results
in \cite{KK14} and \cite{DGM16} have proven that degrees $\bfd$
for which $\bfd'\geq_T\bfz''$ form the spectrum of the real closed field
$\Rz$ containing precisely the computable real numbers.
The article \cite{DGM16} went further, investigating
the archimedean real closed field containing all real numbers computable from degrees
in a given countable ideal $\mathcal{I}$ within the Turing degrees.  The goal
of this article is to investigate the possible spectra of real closed fields further.

Our work here is organized roughly by scope.  After introducing our definitions and conventions,
we begin in Section \ref{sec:arch} with basic results, which apply to archimedean real closed fields.
Theorem \ref{thm:spec} describes conditions for a degree to belong to the
spectrum of such a field.  In certain respects, this theorem can be viewed as dividing
the problem in two:  the degree needs to be able to enumerate the Dedekind cuts
realized in the field, and its jump needs to be able to decide the algebraic dependence
relation on those cuts.  We provide examples to distinguish these two requirements.

In Section \ref{sec:noncone}, we consider real closed subfields
of the field $\Rz$ of all computable real numbers.
We show that it is possible to take any nonzero c.e.\ degree $\bfd$
and give a real closed subfield of $\Rz$ whose spectrum contains $\bfd$
but not $\bfz$.  In particular, this can be done even when $\bfd$ is low,
meaning that the spectrum cannot be defined merely by conditions on jumps of degrees.
The results are reminiscent in certain respects of known theorems about spectra
of linear orders, but some of the results for linear orders remain open (and
appear more challenging) when one asks about real closed subfields of $\Rz$.

In Section \ref{sec:nonarch}, we continue on to nonarchimedean real closed fields.
Here many more questions remain unanswered.  However, in Theorem
\ref{thm:nocompstd} we give the first example (we believe)
of a computable real closed field $F$ whose residue field has no computable presentation.
We also use this field to establish that the question of the spectrum of a real closed
field depends on more than just the spectrum of the residue field and the spectrum
of the derived linear order of the positive infinite multiplicative
classes:  the $F$ built here has the same derived linear order
as another real closed field built earlier in this article, and the residue fields
of the two both have the same spectrum, yet the two fields themselves have distinct spectra.

A theorem of Madison in \cite{M70} shows that for every ordered field $(F,<)$,
the real closure $RC(F)$ of $F$ has a presentation computable from the atomic diagram
of $(F,<)$, and indeed the process of computing the real closure from
this diagram is uniform.  Therefore, $\spec{F,<}\subseteq\spec{RC(F)}$.
However, this containment can be proper.  We leave for another time
the question of degree spectra of ordered fields.

\section{Computable Dedekind Cuts}
\label{sec:Dedekind}

We fix a computable bijection between the sets $\omega$
(of nonnegative integers) and $\Q$ (of rational numbers,
viewed as equivalence classes in $\omega\times(\omega-\{0\})$).
Thus we may speak of computable and c.e.\ subsets of $\Q$,
as well as of $\omega$.

We will need to deal with three distinct kinds of Dedekind cuts.
\begin{defn}
\label{defn:Dedekind}
A \emph{right-leaning Dedekind cut} consists of two subsets $A$ and $B$ of $\Q$ such that:
\begin{itemize}
\item
$A\cup B=\Q$ and $A\neq\emptyset$ and $B\neq\emptyset$; and
\item
for every $q\in A$ and $r\in B$, we have $q<r$; and
\item
$A$ has no greatest element (under the usual order $<$ on $\Q$).
\end{itemize}
The number \emph{realized}, or \emph{defined}, by this cut
is the unique real number $x$ in $\bigcap_{q\in A, r\in B} (q,r]$.
For a \emph{left-leaning Dedekind cut}, we alter the final item
to stipulate that $B$ has no least element, but now $A$ is allowed
to have a greatest element.  In a \emph{strict Dedekind cut},
neither has a least nor a greatest element, and we weaken the first item,
which now requires that $A\neq\emptyset\neq B$ and that $(\Q-(A\cup B))$
be either empty or a singleton.
Of course, the distinction between the three kinds of Dedekind cuts
is trivial when $x$ is irrational.  For rational $x$, the different
kinds simply specify whethere $x$ itself should lie: in $B$, in $A$,
or in neither.

If $A$ and $B$ are both $\bfd$-computably enumerable
(or equivalently, if they are both $\bfd$-computable),
then this is a \emph{$\bfd$-computable Dedekind cut}.  The term
applies to all three kinds.

A \emph{$\bfd$-computable enumeration of Dedekind cuts}
consists of two sequences $\{ A_n\}_{n\in\omega}$ and $\{ B_n\}_{n\in\omega}$,
both uniformly $\bfd$-c.e., such that, for every $n$, $(A_n,B_n)$
forms a Dedekind cut.
\end{defn}

When $A=\cup_s A_s$ and $B=\cup_s B_s$ are $\bfd$-computable
enumerations of these sets, with every $A_s$ and $B_s$ finite and nonempty,
it is often convenient to describe the Dedekind cut $(A,B)$
by a sequence of nested intervals describing the possible
values of the real filling the cut: $(a_s,b_s)$ for a strict cut,
or $(a_s,b_s]$ or $[a_s,b_s)$ if the cut leans right or left,
where $a_s=\max A_s$ and $b_s=\min B_s$.
Of course, many distinct sequences of intervals
will correspond to the same cut $(A,B)$, depending on the
enumerations used for $A$ and $B$.

The following result appears as Lemma 2.1 in \cite{DGM16}.
\begin{lemma}
\label{lemma:Friedberg}
Let $\{ (A_n,B_n)\}_{n\in\omega}$ be a $\bfd$-computable
enumeration of Dedekind cuts.  Then there exists another
$\bfd$-computable enumeration $\{ (C_k,D_k)\}_{k\in\omega}$
realizing exactly the same real numbers, but such that, for all
$j<k$, $(C_j,D_j)$ and $(C_k,D_k)$ realize distinct real numbers.
(We call the latter a \emph{Friedberg enumeration}, in honor
of the originator of the similar theorem for c.e.\ sets.)

Moreover, there is a procedure for building the
Friedberg enumeration from the original enumeration,
uniformly for all enumerations $\{ (A_n,B_n)\}_{n\in\omega}$
and uniformly relative to the degree $\bfd$.  Also, there is a
$\bfd$-computable function $f$ (uniformly, again) such that, for every $k$,
$(C_k,D_k)$ realizes the same real number as $(A_{f(k)},B_{f(k)})$.
\qed\end{lemma}
This theorem applies to all three kinds of Dedekind cuts.
For a Friedberg enumeration, one simply waits to enumerate a cut
until it has distinguished itself from the finitely many
cuts preceding it in the original enumeration.  (This is
far easier than Friedberg's own result for c.e.\ sets.)

Computable Dedekind cuts $(A,B)$ are those for which
both $A$ and $B$ are computable -- or equivalently, both are c.e.
In this case, the unique real number filling the cut $(A,B)$ is said
to be computable.  Computable real numbers form a simple bridge
connecting the computable Dedekind cuts with the computable
subsets of $\omega$.
\begin{lemma}
\label{lemma:allnonstrict}
There is an effective bijection between computable subsets $S\subseteq\omega$
and computable non-strict Dedekind cuts of real numbers in $[0,1]$.
\end{lemma}
``Non-strict'' means that, when $x\in (0,1)$ is rational, we include
both the left-leaning and the right-leaning cut defined by $x$.
To make things perfect, we include the right-leaning cut of $0$
and the left-leaning cut of $1$.  Then the bijection is defined
recursively by letting $n$ lie in $S$ just if
$$\frac1{2^{n+1}}+\sum_{m\in S~\&~m<n}\frac1{2^{m+1}}$$
lies in the left side of the cut.

\begin{lemma}
\label{lemma:bridge}
Let $S\subseteq\omega$.  Then $S$ is computable if and only if the
real number $r_S=\sum_{n\in S} \frac1{2^{n+1}}\in[0,1]$ is computable.
\qed\end{lemma}
It is possible here to have $r_S=r_T$ for distinct sets $S$ and $T$
(e.g., for $S=\{ 1\}$ and $T=\omega-\{0,1\}$),
but only if one set is finite and the other cofinite.  With this lemma,
the next result follows quickly, using a theorem proven by Jockusch
in \cite{J72}.  Again, this result was established in \cite{DGM16},
in a more general form:  instead of considering sets $\leq_T\bfc$,
they considered all sets in a countable ideal $\mathcal{I}$ of Turing degrees.

\begin{prop}[see \cite{J72} and Theorem 1.1 of \cite{DGM16}]
\label{prop:bridge}
Fix any Turing degree $\bfc$.  Then for each Turing degree $\bfd$,
the following are equivalent.
\begin{itemize}
\item
$\bfd$ can enumerate the set of all $\bfc$-computable nonstrict
Dedekind cuts of numbers in $[0,1]$;
\item
$\bfd$ can enumerate the $\bfc$-computable sets;
\item
$\bfd'\geq\bfc''$, i.e., $\bfd$ is \emph{high relative to $\bfc$}.
\qed\end{itemize}
\end{prop}
The equivalence of the last two conditions, which is not trivial,
is Jockusch's result,
relativized to $\bfc$, while the equivalence of the first two follows
from Lemma \ref{lemma:allnonstrict}.

\section{Known Spectra of Real Closed Fields}
\label{sec:known}

The next result was proven (with $\bfc=\bfz$)
by Korovina and Kudinov \cite{KK14}, as well as in \cite{DGM16}.

\begin{thm}
\label{thm:R_c}
For each Turing degree $\bfc$, 
the spectrum of the field $\Rc$ of all $\bfc$-computable real numbers
contains precisely those degrees $\bfd\geq\bfc$ which are high relative
to $\bfc$, i.e., with $\bfd'\geq\bfc''$.  Moreover, these are the only
degrees capable of enumerating the Dedekind cuts realized in $\Rc$.
\qed\end{thm}

Of course, each upper cone of Turing degrees $\geq_T\bfc$
is the spectrum of the field $\Q(x)$,
where $x$ is any real number of degree $\bfc$,
and is also the spectrum of the real closure of $\Q(x)$.
We can combine this with the
preceding idea as follows, to obtain spectra which are neither
upper cones nor jump-preimages of upper cones.
Recall the definition:  an ordered field $F$ is \emph{archimedean} if,
for every element $x\in F$, there is a natural number $n$ with $x<n$.
\begin{thm}
\label{thm:combo}
For every pair of Turing degrees $\bfc_0$ and $\bfc_1$,
there exists an archimedean real closed field $F$ with
$$ \spec{F} = \set{\bfd}{\bfc_0\cup\bfc_1\leq\bfd~\&~\bfc_1''\leq\bfd'}.$$
\end{thm}
\begin{pf}
Fix a real number $x$ of degree $\bfc_0$, and let $F$ be the real closure
of $\R_{\bfc_1}(x)$, where $\R_{\bfc_1}$ contains all $\bfc_1$-computable
real numbers.  Clearly every $\bfd\in\spec{F}$ computes the cut of $x$,
so $\bfd\geq\bfc_0$.  Also, with a $\bfd$-oracle, we may start enumerating
all non-strict Dedekind cuts (both left-leaning and right-leaning) of elements
$y\in F$ with $0\leq y\leq 1$.  Whenever we see that, for some $y\in F$,
$x$ is algebraic in $F$ over $\Q(y)$, we change the enumerations of the two cuts
of $y$ so that they enumerate the cuts of some convenient rational number instead.
Thus, in the end, our enumeration contains exactly the non-strict cuts
in $[0,1]$ realized in the ground field $\R_{\bfc_1}$, i.e., the $\bfc_1$-computable
cuts in this interval.  By Proposition \ref{prop:bridge}, $\bfd$ must be
$\geq\bfc_1$ and high relative to $\bfc_1$ (that is, $\bfd'\geq\bfc_1''$).

Conversely, if $\bfd\geq\bfc_1$ and $\bfd'\geq\bfc_1''$, then $\bfd$ computes
a copy of $\R_{\bfc_1}$, by Theorem \ref{thm:R_c}, and if also $\bfd\geq\bfc_0$,
then we can extend this copy to a presentation of $\R_{\bfc_1}(x)$
as an ordered field.  (If $\bfc_0\leq\bfc_1$, this is just $\R_{\bfc_1}$ itself;
whereas if not, then $x$ must be transcendental over $\R_{\bfc_1}$,
and so the field arithmetic is just that of a purely transcendental
extension of $\R_{\bfc_1}$.)  By Madison's theorem from \cite{M70},
the real closure $F$ of $\R_{\bfc_1}(x)$ is also $\bfd$-computably
presentable.
\qed\end{pf}

So far, therefore, the spectra of archimedean real closed fields
we have met are exactly those named in Theorem \ref{thm:combo}:
the set of degrees in some upper cone satisfying a particular
highness condition.  (The condition $\bfd'\geq\bfc''$ defines highness relative
to $\bfc$; the set of degrees satisfying this condition may also be described
as the preimage of the upper cone above $\bfc''$ under the jump operation.)
Adjoining finitely many elements $x_1,\ldots,x_n$ of incomparable degree to a given
real closed field does not create any further spectra: the upper cone
now is simply that above the join of the degrees of $x_1,\ldots,x_n$.
However, in Section \ref{sec:noncone} we will produce
archimedean real closed fields with spectra distinct from those
described in Theorem \ref{thm:combo}.

\section{Archimedean Real Closed Fields}
\label{sec:arch}

Proposition \ref{prop:bridge} shows that for the countable
real closed fields $\Rc$, the ability to enumerate the cuts of
the real numbers in $\Rc$ is equivalent to the ability to
present the field $\Rc$.  In general, however,
having a computable Friedberg enumeration
of the Dedekind cuts realized in an archimedean real closed field $F$
is not \emph{a priori} sufficient to yield a computable presentation
of $F$.  In fact, it is not difficult to give a computable Friedberg
enumeration $(~\la a_{i,s},b_{i,s}\ra~)_{i,s\in\omega}$ of the cuts
realized in the real closure of $\Q$ in such a way that the addition
function (mapping a pair $(i,j)$ to the unique $k$ such that
$\lim_s a_{k,s}=\lim_s (a_{i,s}+a_{j,s})$) and the multiplication function
each have degree $\bfz'$.  ($\bfz'$ is readily seen to be sufficient
to compute these functions, given any computable Friedberg enumeration
for an archimedean real closed field $F$ as described above.)
In Theorem \ref{thm:Sigma2} we will give a separate example.

The basic criteria for belonging to the spectrum of a real closed field
are given by the following theorem.

\begin{thm}
\label{thm:spec}
Let $F$ be an archimedean real closed field.  Then, for each
Turing degree $\bfd$, the following are equivalent.
\begin{enumerate}
\item
$\bfd\in\spec{F}$.
\item
There is a $\bfd$-computable Friedberg enumeration of the strict Dedekind cuts realized in $F$,
in which the algebraic dependence relation on those cuts is $\bfd'$-decidable.
\item
There is a $\bfd$-computable Friedberg enumeration of a set of Dedekind cuts
such that the real numbers filling those cuts form a transcendence basis for $F$.
\end{enumerate}
\end{thm}
(In (2), one can substitute either right-leaning or left-leaning Dedekind cuts
for strict ones, yielding two more equivalent conditions.)

\begin{pf}
It is immediate that $(1)\implies (2)$.  A $\bfd$-computable presentation of $F$
gives a Friedberg enumeration of the cuts realized in $F$, just by enumerating the elements
of the presentation and listing the cuts in which they lie.  Since $\bfd$ also computes
the field operations on these cuts, the dependence relation on these cuts
these elements is $\Sigma^{\bfd}_1$, hence $\bfd'$-decidable, proving (2).  

It is also clear that $(1)\implies (3)$ when $F$ has finite transcendence degree
over $\Q$:  just enumerate the cuts for any transcendence basis.
For a field $F$ of infinite transcendence degree over $\Q$,
start with a $\bfd$-computable presentation of $F$.  We define
cuts $C_0,C_1,\ldots$.  At stage $0$, we begin enumerating $C_0$
as the cut of the element $x_0$ in the domain $\{ x_0,x_1,\ldots\}$ of $F$,
saying that $C_0$ is being \emph{guided} by $x_0$.
At stage $1$ we start enumerating $C_1$ as the cut of $x_1$,
and so on.  However, we also search at the same time for algebraic
relations in $F$ among these elements.  If we ever see that $x_0$
is algebraic over $\Q$, then we drop $x_0$ and perform a ``left-shift''
as follows.  For the element $x_{i_1}$ currently guiding $C_1$,
find a rational $q_1$ such that $x_{i_1}+q_1$ lies in $C_0$ as currently
defined, and let $x_{i_1}+q_1$ guide $C_0$ from now on.
Likewise, for the $x_{i_2}$ currently guiding $C_2$, find $q_2\in\Q$
with $x_{i_2}+q_2$ in the current interval defined by $C_1$,
and let $x_{i_2}+q_2$ guide $C_1$ from now on, and so on.
For the greatest $C_s$ currently being enumerated, take the next available
$x_i$ and start using it (plus a rational) to guide $C_s$ from now on.

Likewise, whenever the element $x_{i_j}$ currently guiding $C_j$ is
found to be algebraic over $\Q(x_0,\ldots,x_{i_j-1})$, we perform a
left-shift for all the cuts $C_j,C_{j+1},\ldots$, without disturbing
$C_0,\ldots,C_{j-1}$.  Since $F$ has infinite transcendence degree,
this ensures that $C_0$ will eventually be the cut of an element
which differs by a rational from the first $x_{i_0}$ transcendental over $\Q$.
Likewise, $C_1$ will be guided by the first $x_{i_1}$ transcendental over
$\Q(x_{i_0})$, and so on for each $C_j$.  This proves (3).

To prove $(3)\implies (1)$, we use the $\bfd$-computable enumeration
of the cuts for a transcendence basis for $F$ to build a $\bfd$-computable
ordered field $\Q(x_0,x_1,\ldots)$, with one $x_i$ for each cut in the enumeration.
The field arithmetic is that of a purely transcendental extension, and the order
is decidable from the enumeration of the cuts.  Then Madison's theorem from \cite{M70}
yields a $\bfd$-computable presentation of the real closure of this field, namely $F$.

We show $(2)\implies(3)$ by the same basic process as 
for $(1)\implies (3)$.  The $\bfd'$-decidability of algebraic dependence on the
$\bfd$-computable enumeration of cuts allows $\bfd$ to approximate the
characteristic function of a transcendence basis.  Sometimes a right-shift
is necessary, by analogy to the left-shifts described above.  If the approximation
says that the first cut is algebraic, but then changes its mind, then $C_0$
must start its enumeration of the first cut (up to a rational difference) all over
again, moving all other cuts one step to the right.  However, the process still
succeeds, yielding a $\bfd$-computable enumeration of the cuts in a
transcendence basis for $F$.
\qed\end{pf}

Theorem \ref{thm:spec} provides a degree-theoretic criterion sufficient
for membership in $\spec{F}$, though not always necessary.

\begin{cor}
\label{cor:necessary}
Let $F$ be an $\bfa$-computable archimedean real closed field,
such that every $x\in F$ is $\bfc$-computable. Then every degree $\bfd\geq\bfc$
with $\bfd'\geq\bfa'\cup\bfc''$ lies in $\spec{F}$.
\end{cor}
\begin{pf}
Given any natural number $x$ (which we regard as an element
of the $\bfa$-computable presentation $F$) and any $i,j\in\omega$,
a $\bfd'$-oracle can decide whether
\begin{align*}
& (W^C_i,W^C_j)\text{~is a strict Dedekind cut in $\Q$}\\
\&~&(\forall q\in\Q)(\forall s)[(q\in W^C_{i,s}\to q<x)~\&~(q\in W^C_{j,s}\to x < q)]\\
\&~&(x\text{~is transcendental over $\Q(0,\ldots,x-1)$ in $F$.})
\end{align*}
(The first condition is $\Pi_2^{\bfc}$.
The second is $\Pi_1^{\bfa}$, using the $<$ relation in $F$,
and so is the third.)
So we have a $\bfd$-computable approximation
to this property, uniformly for all $x$, $i$, and $j$.
Of course, for each $i$ and $j$, at most one $x$
can make the property true.  So we use a $\bfd$-computable
finite injury procedure, building a cut in $\Q$ to approximate
$(W_i^C, W_j^C)$ as long as it appears that some specific $x$
makes the property true for this (potential) cut, and
doing a left-shift or a right-shift, as in Theorem \ref{thm:spec},
if the approximation ever changes its mind for this particular $(i,j,x)$.
This yields a $\bfd$-computable list of all the cuts in
a transcendence basis for $F$, so $\bfd\in\spec{F}$
by Theorem \ref{thm:spec}.
\qed\end{pf}

Between the two criteria given by this theorem for membership in $\spec{F}$,
condition (3) seems by far the cleaner and more useful.
We include condition (2) because it suggests the dual
requirements in presenting $F$:  enumerating the cuts, and computing
the field operations.  We now give an example which, in concert with
Theorem \ref{thm:R_c}, shows how each of these dual requirements
can take precedence.  Both examples have the same spectrum.  However,
the new example is a real closed field $R_S$ which satisfies the first
part of (2) with $\bfd=\bfz$, but not the second.  Roughly speaking, in this field,
it is the field arithmetic which prohibits the existence of a computable copy.
In Theorem \ref{thm:R_c}, in contrast, with the field
$\Rc$ of $\bfc$-computable real numbers, the field
arithmetic causes no problems:  the ability to give a Friedberg enumeration
of the cuts realized in $\Rc$ was sufficient for a degree to lie in its spectrum.

\begin{thm}
\label{thm:Sigma2}
For every $S\in\Sigma^0_2$, there is a countable archimedean real closed field
$R_S$ with a computable Friedberg enumeration of the cuts realized in $R_S$,
such that $\spec{R_S}=\set{\deg{D}}{S\leq_T D'}$.
\end{thm}
When $S$ has degree $\bfz''$, this spectrum has already been realized
by the field $\R_{\bfz}$ of all computable real numbers (see Theorem \ref{thm:R_c}).
The point of this theorem is that, unlike $\R_{\bfz}$, $R_S$ has a
computable enumeration of its cuts; the complexity of $\spec{R_S}$
stems entirely from the difficulty of the field operations in $R_S$.
\begin{pf}
Fix a computable $1$-reduction $f$ from $S$ to the $\Sigma^0_2$-complete
set $\Fin$.  First we will uniformly
enumerate Dedekind cuts $\set{ (a_{e,s},b_{e,s})}{i\in\omega}$
such that, for each $e$, $a_e=\lim_s a_{e,s}$ is algebraic over $\Q$ if and only if
$e\in S$ (that is, just if $W_{f(e)}$ is finite).  Moreover, if $e\notin S$, then
$a_e$ will be transcendental over the subfield $\Q(a_0,\ldots,a_{e-1})$.

Let $p_0,p_1,\ldots$ list all nonzero polynomials in $\Z[X]$.  At stage $0$
of the construction, we choose each interval $I_{e,0}=(a_{e,0},b_{e,0})=(e,e+1)$,
and fix the least index $n_{e,0}$ such that $p_{n_{e,0}}(X)$
has a root in $I_{e,0}$.

At stage $s+1$, for each $e$, we check whether $W_{f(e),s+1}=W_{f(e),s}$.
If so, then we take $I_{e,s+1}$ to be a subinterval of $I_{e,s}$ as follows. Suppose $I_{e,s}=(q_{e,s}, q'_{e,s})$, with $q_{e,s}, q'_{e,s}\in \Q$. Define $I_{e,s+1}=(q_{e,s+1}, q'_{e,s+1})$, with $q_{e,s}, q'_{e,s}\in\Q$ satisfying $q_{e,s} < q_{e,s+1}< q'_{e,s+1} <q'_{e,s}$ and $(\exists a)$ $(q_{e,s+1}< a < q'_{e,s} \text{ } \&~p_{n_{e,s}}(a)=0)$. We keep $n_{e,s+1}=n_{e,s}$.  

If $W_{f(e),s+1}\neq W_{f(e),s}$, then suppose $I_{e,s}=(q_{e,s}, q'_{e,s})$,
with $q_{e,s}, q'_{e,s}\in \Q$. Define $I_{e,s+1}=(q_{e,s+1}, q'_{e,s+1})$,
with $q_{e,s}, q'_{e,s}\in\Q$ satisfying
\begin{itemize}
\item
$q_{e,s} < q_{e,s+1}< q'_{e,s+1} <q'_{e,s}$,
\item
$(\forall a)$ $(q_{e,s+1}< a < q'_{e,s} \rightarrow p_{n_{e,s}}(a)\neq 0)$, and
\item
for each of the first $s$ nonzero polynomials $r\in\Z[X_0,\ldots,X_{e-1},X]$
and each $x_0\in (a_{0,s},b_{0,s}),\ldots,x_{e-1}\in (a_{e-1,s},b_{e-1,s})$
no root of $r(x_0,\ldots,x_{e-1},X)$ lies in $I_{e,s+1}$.
\end{itemize}
(This might not be possible
at this stage, but it is decidable whether it is possible or not.
If for some $r$ it is impossible, then we ignore that $r$ at this stage.
Eventually the preceding intervals $(a_{i,s},b_{i,s})$ will contract
enough that it will be possible.)
We define $n_{e,s+1}$ to be the least index such that the polynomial
$p_{n_{e,s+1}}(X)$ does have a root in $I_{e,s+1}$.

This defines the entire sequence of cuts $\set{(a_{i,s},b_{i,s})}{i\in\omega}$,
all with distinct limits. Let $R_S$ be the real closed field generated by the elements $a_e$
realized by these cuts.  We can extend this sequence to an effective
Friedberg enumeration of all cuts realized in $R_S$.  Let $R$ be the
algebraic dependence relation on the cuts in this enumeration.

It is clear from this construction that $f(e)\in\Fin$ if and only if the
element $a_e=\lim_s a_{e,s}$ of $R_S$ is algebraic over $\Q$.
Thus $S\leq_1 R$.  Conversely, since every $e\notin\Fin$
has a cut $I_e=\lim_s I_{e,s}$ trancendental over all the cuts
$I_j$ with $j<e$, we can decide $R$ from an $S$-oracle;
hence $R\equiv_T S$.

But now, suppose that $E$ is a copy of $R_S$ of Turing degree $\bfd$.
Then, with a $\bfd'$-oracle, we can take an arbitrary $e$ and find the
element $z_e$ of $E$ which realizes the same cut in $E$ that $a_e$
realizes: the cuts in $E$ are $\bfd$-computably enumerable, and
the enumeration of the cuts $I_e$ was $\bfz$-computable,
so a $(\bfd\cup\bfz)'$-oracle is all that is needed.
Then we determine (still using $\bfd'$) whether $z_e$ is transcendental
over $\Q$ in $E$.  Therefore, $\bfd'$ computes $S$.

Conversely, if $\bfd'$ computes $S$, then $\bfd$ satisfies
item (2) of Theorem \ref{thm:spec}, using our computable
enumeration of the cuts in $R_S$, so $\bfd\in\spec{R_S}$.
\qed\end{pf}

In addition to showing the relevance of the field operations
to the spectrum, Theorem \ref{thm:Sigma2} establishes that
the jump-preimage $\set{\bfd}{\bfd'\geq\deg{S}}$ can be
a spectrum even when $\deg{S}$ is not a second jump, i.e., when
$S\not\geq_T\emptyset''$.  The fields $\Rc$,
with their spectra $\set{\bfd}{\bfd'\geq_T\bfc''}$, did not accomplish this.
\begin{cor}
\label{cor:noncomputable}
For every nonlow $\Delta^0_2$ set $U$, there exists a
$U$-computable archi\-medean real closed subfield of $\Rz$
which has no computable presentation.
\end{cor}
\begin{pf}
Apply Theorem \ref{thm:Sigma2} to $U'$.
\qed\end{pf}

\section{New Spectra of Real Closed Fields}
\label{sec:noncone}

It is clear that, if an archimedean real closed field $F$ contains
an element whose Dedekind cut has Turing degree $\bfc$,
then the spectrum of $F$ is contained within the upper cone above $\bfc$.
So far, the converse has also been true:  among the archimedean
real closed fields $F$ that we have seen so far,
each one whose spectrum lies within the upper cone above any degree $\bfc$
contains an element whose Dedekind cut has that degree.
This converse does not hold in general, but
a theorem of Knight \cite[Thm.\ 1.4]{K86} shows that, whenever
the spectrum of a real closed field $F$ lies in the upper cone above
$\deg{C}$, there exists some finite tuple $\avec$ from $F$ such that
both $C$ and its complement $\Cbar$ are $e$-reducible to the existential
theory of $(F,\avec)$.  (This result was proven independently by Soskov,
in addition to the proof by Knight in \cite{K86}.)  Since the theory of
real closed fields is decidable, the $\Sigma_1$-theory of $(F,\avec)$
is enumerable from the join of the Dedekind cuts of the elements $\avec$,
and so $C$ is computable in this join.  We record this property here.

\begin{prop}[following Knight and Soskov]
\label{prop:meet}
Let $F$ be an archimedean real closed field.  Then for every set $C$
with $\spec{F}\subseteq\set{\bfd}{\deg{C}\leq_T\bfd}$, there exists a 
finite tuple $\avec$ of elements of $F$ such that the join of the Dedekind cuts
of the elements $\avec$ computes $C$.
\qed\end{prop}

The fact that no fixed finite number of elements suffices is shown by a finitary
forcing argument: one builds the real closure of a field $\Q(x_1,\ldots,x_n)$,
using Friedberg-Muchnik requirements to ensure that these $x_i$ have
cuts of pairwise-incomparable degree, mixed with requirements
that, if the cut of an element $y$ computes the cut of any $x_i$,
then $y$ lies in the real closure of $\Q(x_i)$.
(If $y$ is not algebraic over $\Q(x_i)$, then we can adjust
the cut of $y$ by adjusting that of some other $x_j$,
so as to satisfy this requirement.)

Separately one can also build the real closure
$F$ of a field $\Q(x_0,x_1,\ldots)$ such that the cut of $x_{i+1}$
is never computable from the join of the cuts of $x_0,\ldots, x_i$,
and thereby show that there need not be any greatest degree $\bfc$
among those whose upper cone contains $\spec{F}$.

So far a stronger condition than Proposition \ref{prop:meet} has also held
of our examples.  Each $F$ so far has the property of \emph{first-jump
equivalence} among those degrees computing all cuts realized in $F$:
if $\bfd_0$ and $\bfd_1$ both compute all these cuts and
$\bfd_0'=\bfd_1'$, then $\bfd_0\in\spec{F}$ if and only
$\bfd_1\in\spec{F}$.
Now we show that more spectra than just these are possible.
By working with subfields $F$ of $\Rz$, we ensure that the requirement
of computing all the cuts in $F$ is trivial.

\begin{thm}
\label{thm:low}
For every c.e.\ set $L>_T\emptyset$,
there exists an $L$-computable real closed subfield $F$
of $\Rz$ with no computable copy.  In particular, this holds
when $L$ is low and noncomputable, in which case
first-jump equivalence fails for $\spec{F}$.
\end{thm}
\begin{pf}
%
Fix computable enumerations of $L=\bigcup_s L_s$, of $\Q=\{q_0,q_1,\ldots\}$,
and of $\Z[X_0,X_1,\ldots]=\{0\}\cup\{ p_0,p_1,\ldots\}$, allowing
nonzero polynomials with arbitrarily many variables among our $p_i$.
Also fix the usual computable numbering of $\omega^2$ and let
$$ D_{\la j,k\ra} = (W_j,W_k)$$
be the $\la j,k\ra$-th pair of c.e.\ subsets of $\Q$.

Clearly it is a $\Pi^0_2$ property for $i\in\Ded$, that is,
for $D_i$ to be a Dedekind cut.  Indeed, it is $\Pi^0_2$-complete:
for a $1$-reduction
from $\Inf$ to $\Ded$, just map each index $e$ to a pair $\la g(e),k\ra$
such that, whenever $W_e$ receives a new element, $W_{g(e)}$ enumerates
the next available rational number $<\sqrt2$, while $W_k$ contains exactly
those rationals $>\sqrt2$.
(It makes no difference here
whether we consider right-leaning, left-leaning, or strict cuts.)
The property $D_{i_0}=D_{i_1}$ is also $\Pi^0_2$-complete.
Each computable sequence $\calS$ of Dedekind cuts can now be
given as $\la D_{f(n)}\ra_{n\in\omega}$ by a computable function $f$.
We will enumerate these sequences effectively as
$\calS_e=\la D_{\phi_e(n)}~:~(\forall m\leq n)~\phi_e(m)\converges\ra$,
Of course, $\calS_e$ may be a finite sequence of cuts, but the set
$$\Ss=\set{e}{e\in\Tot~\&~
\calS_e\text{~consists of algebraically independent Dedekind cuts}}$$
is still $\Pi^0_2$.  By Theorem \ref{thm:spec}, we may think of $\Ss$
as a listing of all computably presentable real closed subfields
of $\R_0$ of infinite transcendence degree, and use a ``chip function''
such that $e\in\Ss$ if and only if $e$ receives infinitely many chips.
We write $\Ss=\{e_0<e_1<e_2<\cdots\}$, knowing this listing to be noneffective.

We now approximate, computably, a specific sequence of Dedekind cuts
$C_m$ by recursion on $m$, writing $x_m$ for the real number
filling the cut $C_m$.  
The goal is to satisfy two types of requirements, for all $m$ and $i$.  
\begin{align*}
\mathcal N_{\la m,i\ra}:~~&\text{If $y_1,\ldots,y_n$ fill the first $n$ cuts of~}\calS_{e_m},~
p_i(y_1,\ldots,y_n,x_m)\neq 0.\\
\P_i:~~&p_i(x_0,\ldots,x_n)\neq 0.
\end{align*}
(In each case, $n$ is determined simply by the number of variables in $p_i$.)
The $\P$-requirements will make the set $\set{x_m}{m\in\omega}$
algebraically independent, while the $\mathcal N$-requirements
collectively will show that each $x_m$ is transcendental
over the set of real numbers realizing cuts in $\calS_{e_m}$.
Satisfying all these requirements therefore will mean that
$\set{x_m}{m\in\omega}$ is a transcendence basis for a real closed field
with no computable presentation, as no $\calS_{e_m}$ is a transcendence
basis for the same field.  We will then build an $L$-computable
enumeration of these cuts and appeal to Theorem \ref{thm:spec}
to prove Theorem \ref{thm:low}.

Since the field $\Rz$ has no computable presentation, Theorems
\ref{thm:spec} and \ref{thm:R_c} show that no sequence $\calS_e$
can give a basis for $\Rz$, and indeed $\Rz$ must have infinite
transcendence degree over each $\calS_{e_m}$.  Therefore, if $\phi_e$ is
total and has image $\subseteq\Ded$,
then there must be some computable real number $x$ independent
from $\{ x_0,\ldots,x_{e-1}\}$ which is filled by no cut in $\calS_e$,
so it is reasonable to hope to satisfy these requirements.  To do so,
we will use the following simple lemma, saying that as the cuts
for $x_1,\ldots,x_n$ close in on their values, we will eventually be able
to define the cut of $x_0$.
\begin{lemma}
\label{lemma:density}
Let $p\in\Z[X_0,\ldots,X_{m}]$ be a nonzero polynomial,
and fix 
an algebraically independent set $\{x_0,\ldots,x_{m-1}\}$.
Then for every $a_{m}<b_{m}$ in $\Q$, there exist rational numbers $a_i<x_i$
and $b_i>x_i$ (for all $i < m$) and $a<b\in\Q$ with $(a,b)\subseteq (a_m,b_{m})$
such that the image of
$$(a_0,b_0)\times\cdots\times (a_{m-1},b_{m-1})\times (a,b)$$
under $p$ does not contain $0$.
\end{lemma}
\begin{pf}
If not, then by continuity of $p$ we would have $p(x_0,x_1,\ldots,x_m)=0$
for all $x_m\in (a_m,b_m)$, so $p(x_0,\ldots,x_{m-1},X_m)$ would be the zero polynomial.
But, writing $p=\sum_i X_m^i\cdot p_i(X_0,\ldots,X_{m-1})$,
we would then have $p_i(x_0,\ldots,x_{m-1})=0$ for every $i$.
With $\{x_0,\ldots,x_{m-1}\}$ independent, this forces $p=0$.
\qed\end{pf}

We now give a computable procedure to approximate the cuts $C_n$
which satisfy our requirements.  This uses our uniform $\Pi^0_2$-guessing procedure
for $\Ss$, along with permitting below the set $L$.
The construction takes place on the tree $T=\omega^{<\omega}$,
and the true path will be the function $m\mapsto e_m$.
We begin by setting $a_{m,0}=0$ and $b_{m,0}=1$ for all $m$,
so that $A_{m,0}=(-\infty,0]$ and $B_{m,0}=[1,+\infty)$,
and $C_{m,0}=(A_{m,0},B_{m,0})$ is our first approximation
to the cut $C_m$.
We adopt the convention that every requirement is
both eligible and active at stage $0$.

At stage $s+1$, we have at most $s+1$ substages, starting with the
substage for the root node $\lambda$, which is eligible at every stage.
The procedure for an eligible node $\sigma_{m,s+1}=(e_{0,s+1},\ldots,e_{m,s+1})$
with $m\leq s$ at stage $s+1$ is as follows.  Write $e=e_{m,s+1}$,
and let $s'$ be the last stage at which either this node or any node
to its left was active.
Consider the requirements
$$ \mathcal N_{\la m,0\ra} \prec \P_{0} \prec\mathcal N_{\la m,1\ra}\prec \P_{1}\prec\cdots .$$
A requirement $\mathcal N_{\la m,i\ra}$ on this list is \emph{currently satisfied} if
either $\phi_{e,s}(n)\diverges$ for some $n<m$, or else (writing $\la j_n,k_n\ra=\phi_e(n)$)
the image of
$$(\max(W_{j_0,s}),\min(W_{k_0,s}))\times\cdots\times
(\max(W_{j_{m-1},s}),\min(W_{k_{m-1},s}))
\times(a_{m,s},b_{m,s})$$
under $p_i$ does not contain $0$.
(These intervals are the ones defined so far by $\calS_e=\calS_{e_{m,s+1}}$, the 
computable sequence of cuts over which $x_m$ is supposed to be made transcendental.)
Likewise, $\P_i$ is \emph{currently satisfied} if either $p_i\notin\Z[X_0,\ldots,X_m]$
or $p_i\in\Z[X_0,\ldots,X_{m-1}]$ or the image of
$$ (a_{0,s+1},b_{0,s+1})\times\cdots\times(a_{m-1,s+1},b_{m-1,s+1})\times
(a_{m,s'},b_{m,s'})$$
under $p_i$ does not contain $0$.

If any of the first $s$ requirements on this list is not currently satisfied,
then we attempt to satisfy the least such $\mathcal{N}_{\langle m,i\rangle}$
or $\mathcal{P}_i$, by searching for rational numbers $a<b$ among $\{ q_0,\ldots,q_s\}$ such that:
\begin{itemize}
\item
$a_{m,s'}\leq a$ and $b\leq b_{m,s'}$; and
\item
for the smallest $l\in L_{s+1}-L_{s'}$, and for all $q_i$ with $i<l$,
$$ (q_i < a_{m,s}\implies q_i < a)~~\&~~(b_{m,s} < q_i\implies b<q_i)$$
(so every $q_i$ which is moved by the $\mathcal N_{\la m,i\ra}$-action below is permitted
by $L$ to be moved); and
\item
we satisfy the requirement in question:
either the image of $p_i$ on
$$(\max(W_{j_0,s}),\min(W_{k_0,s}))\times\cdots\times
(\max(W_{j_{m-1},s}),\min(W_{k_{m-1},s})) \times(a,b)$$
does not contain $0$ (so we satisfy $\mathcal N_{\la m,i\ra}$); or else 
the image of $p_i$ on
$$ (a_{0,s+1},b_{0,s+1})\times\cdots\times(a_{m-1,s+1},b_{m-1,s+1})\times(a,b)$$
under $p_i$ does not contain $0$ (so we we satisfy $\P_i$).
\end{itemize}
For the first such requirement and the least witnesses $\la a,b\ra$,
we define $a_{m,s+1}=a$ and $b_{m,s+1}=b$.  In this case,
$\sigma_{m,s+1}$ is \emph{active} at this stage, on behalf of $e$:
assuming $m<s$, we make the node $\sigma_{m,s+1}\widehat{~}e$ 
eligible at the next substage, where $e$ is
the least number $>e_{m-1,s+1}$ which has received a chip
since the last stage at which $\sigma_{m,s+1}$ was active.  
In case either $a_{m,s+1}<a_{m,s}$
or $b_{m,s+1}>b_{m,s}$, then we say that $C_m$ was \emph{redefined}
at this stage; otherwise $C_m$ was only \emph{refined} here, possibly trivially.
If no such $a$ and $b$ exist for any of the first $s$ requirements, then we keep
$a_{m,s+1}=a_{m,s'}$ and $b_{m,s+1}=b_{m,s'}$ and end the stage
right here.  In this case $\sigma_{m,s+1}$ was only \emph{eligible} at this stage,
not active, and none of its successors was eligible.
Finally, if $m=s$, then no successor is eligible and the stage ends here,
even if $\sigma_{m,s+1}$ was active.

This completes the construction, and we define $F$ to be the real closed
field generated by the real numbers $x_m$ filling the cuts $C_m$ for each $m$.
(Of course, it remains to prove that the $C_m$ really are Dedekind cuts,
and are computable.)  Since the approximations to these cuts were
redefined only when $L$ permitted such redefinition, the usual
permitting argument shows that $F$ will be $L$-computable,
once we have seen it to be a real closed field.

\begin{lemma}
\label{lemma:satisfaction}
For every $m$, the node $\sigma_m=(e_0,\ldots,e_m)$
(from the listing $\Ss=\{e_0<e_1<e_2<\cdots\}$ defined above)
is the leftmost node at level $m$ to be active at infinitely many stages,
and $C_m$ is the Dedekind
cut of a real number $x_m$ which is transcendental both over $\{x_0,\ldots,x_{m-1}\}$
and over all the real numbers realized by cuts in $\calS_{e_m}$.

\end{lemma}
\begin{pf}
By induction on $m$, we may assume that there are infinitely many stages
at which $\sigma_{m-1}$ is active.  By our chip procedure, there must be a stage
$s_0$ such that no node $\sigma$ to the left of $\sigma_m$ is eligible
after stage $s_0$.  However, $e_m$ itself receives infinitely many chips,
so $\sigma_m$ must be eligible at infinitely many stages.
\comment{To see that
$\sigma_m$ is active at infinitely many of these stages, notice that
whenever $C_{m,s}$ has end points $a_{m,s}, b_{m,s}$ at one of these
stages, each $q_j\in (a_{m,s}, b_{m,s})$ gives a polynomial $p_i=X_m-q_j$
for which we will eventually be able to satisfy requirement $\P_i$
simply by refining $C_m$ (without redefining it), and so, regardless
of $L$-permission, $\sigma_m$ will certainly be active at infinitely many stages.}

Now we claim that every requirement in the sequence
$$ \mathcal N_{\la m,0\ra} \prec \P_{0} \prec\mathcal N_{\la m,1\ra}\prec \P_{1}\prec\cdots $$
will eventually be satisfied.  (In the language above, it will be \emph{currently satisfied}
at cofinitely many stages.)  If not, then there is a least requirement $\mathcal{R}$
for which it fails.
Moreover, in this case $\sigma_m$ never acts at any stage after all of the
higher-priority requirements than $\mathcal{R}$ are satisfied,
since thereafter the construction will always identify $\mathcal{R}$ as the next
one needing satisfaction and will not consider anything of lower priority.
In particular, at all stages when $\sigma_m$ is eligible, the value $s'$
in the construction will be the same.
$C_m$ will continue to be refined and even redefined at stages
when nodes to the right of $\sigma_m$
act, but none of those nodes is allowed to move the interval of $C_m$
outside the interval chosen by $\sigma_m$.  So, from the point of view
of $\sigma_m$, no subsequent redefinition of the cut takes place.

Suppose $\P_i$ is this least requirement $\mathcal{R}$.
Now $(C_0,\ldots,C_{m-1})$ must be a sequence of algebraically
independent Dedekind cuts in $\Ss$, by inductive hypothesis,
and so we do eventually reach stages with all of 
$ (a_{0,s+1},b_{0,s+1})\times\cdots\times(a_{m-1,s+1},b_{m-1,s+1})$ defined.
Moreover, each of these intervals must continue to shrink down to radius $0$
as we go through the construction.  Lemma \ref{lemma:density}
shows that eventually an $(a,b)$ must appear which would satisfy $\P_i$
and will also satisfy $(a,b)\subseteq (a_{m,s'},b_{m,s'})$ for the stage $s'$
(which is fixed for all these $s$, as noted above).  Moreover,
every subinterval $(a',b')\subseteq (a,b)$ would have these same properties.
However, the construction never acts to use any of these subintervals,
so each of them must violate the $L$-permitting condition.  We claim that this
gives us a way to compute $L$.  We can be sure that, whenever
we see a new subinterval $(a,b)$ of $(a_{m,s'},b_{m,s'})$ appear which
could be used to satisfy $\P_i$ at a stage $s$, $L$ will not permit this satisfaction,
and so $L_s\res j=L\res j$, where $j$ is least such that $q_j\in (a,b)$.
Since infinitely many $j$ lie in $(a_{m,s'},b_{m,s'})$, and since each of them
has $j$ minimal within some small $(a,b)$, this yields a method of
deciding membership in $L$, contrary to our hypothesis that $L>_T\emptyset$.
Therefore, eventually one of these subintervals $(a,b)$ must be acted
on by $\sigma_m$, and thereafter $\P_i$ will always be satisfied, since
$\sigma_m$ never acts again on behaf of any higher-priority requirement.

An argument in the exact same style applies if some requirement
$\mathcal N_{\la m,i\ra}$ is this least $\mathcal{R}$.  Now the intervals $(a_{i,s+1},b_{i,s+1})$
are replaced by intervals $(\max(W_{j,s}),\min(W_{k,s}))$ given by $\calS_{e_m}$,
but since $e_m\in\Ss$, we know that these intervals must also form Dedekind cuts
$(W_j,W_k)$, and the rest of the argument is exactly the same.
It now follows that $C_m$ satisfies every requirement on our list.
This proves that $C_m$ really is a Dedekind cut, as follows.
For every rational $q$,
the requirement $\P_i$ corresponding to the polynomial $p_i=X_m-q$
is satisfied, putting $q$ on one side or the other of $C_m$.  
(This shows that the intervals $a_{m,s},b_{m,s})$ really do shrink to $0$
as $s\to\infty$, as was mentioned above.)

These $\P_i$
collectively also show that $\sigma_m$ acts at infinitely many stages.
In turn, the $L$-permitting now ensures the convergence of our Dedekind cuts
$C_m$:  no rational $q_j$ can move out of either side of the cut $C_m$
at any stage after the first stage $s$ at which $\sigma_m$ is active
and $L_s\res j=L\res j$.  Thus our induction is finished.
\qed\end{pf}
Lemma \ref{lemma:satisfaction} is all we need to see that
the real closure $F$ of the ordered field $\Q(x_0,x_1,\ldots)$
is $L$-computable.  Since each $x_m$ is transcendental over
$\calS_{e_m}$, no computable real closed field can be isomorphic
to $F$, by Theorem \ref{thm:spec}.  It remains to show that
every cut $C_m$ is computable, so that $F\subseteq\Rz$.

To compute $C_m$, we need two finite pieces of information.  These are the
string $\sigma_m$ (equivalently, the first $(m+1)$ elements
of $\Ss$) and the least stage $s_0$ after which no $\sigma$
to the left of $\sigma_m$ ever acts again.  Since these data
cannot be computed uniformly in $m$, our proof here does not
show the field $F$ to be computable; we are only proving that
each individual cut $C_m$ is computable (and therefore that all
cuts realized in $F$ are computable, since each of them is the cut
of a real number algebraic over a finite subset of $\{ x_0,x_1,\ldots\}$).

Given $\sigma_m$ and $s_0$, we run the procedure above for building
all the cuts $C_0,C_1,\ldots$.  (This procedure is effective, of course,
although it does not build computable enumerations of the cuts, due to
its occasional redefinitions of various cuts.  It is also uniform in the
enumeration of $L$.)  Once we have passed stage $s_0$, we watch for stages
$s$ at which $\sigma_m$ is active.  Lemma \ref{lemma:satisfaction}
showed that there are infinitely many of these stages, and at each one,
the current values $(a_{m,s},b_{m,s})$ of the end points of the interval $C_{m,s}$
are correct:  all subsequent actions, either by $\sigma_m$ or by nodes
to its right, respect the end points
from the previous stage at which $\sigma_m$ acted, and no $\sigma$
to the left of $\sigma_m$ will ever interfere, since they never act again.
Thus, we have an effective enumeration of each side of $C_m$,
and so indeed $F$ is a subfield of $\Rz$ as required.
\qed\end{pf}

The proof of Theorem \ref{thm:low} relativizes to an arbitrary degree $\bfc$,
yielding an immediate generalization.
\begin{cor}
\label{cor:low}
For every pair of Turing degrees $\bfc<_T\bfd$ with $\bfd$ c.e.\ in $\bfc$,
there is a real closed subfield $F\subseteq\Rc$ whose spectrum contains
$\bfd$ but not $\bfc$.  (We can also ensure that $F$ contains an element
whose Dedekind cut has degree $\bfc$.)
\end{cor}
\begin{pf}
The relativization to $\bfc$ works with no trouble.  To ensure the last condition,
we begin the construction by adjoining to $F$ a single real number of degree $\bfc$:
fix an appropriate $i$ and $j$ with $W_i^C\equiv_T C$, and include a real with cut
$(W_i^C,W_j^C)$, with higher priority than any requirement.
\qed\end{pf}

Theorem \ref{thm:low} leaves the situation for real closed subfields of $\Rz$
essentially where linear orders were after the work of Jockusch and Soare in \cite{JS91}:
each individual c.e.\ degree except $\bfz$ lies in the spectrum of some such structure
with no computable copy.  For linear orders, Downey and Seetapun independently
extended this result (in unpublished work) to all nonzero $\Delta^0_2$-degrees,
and then a result in \cite{M01} showed that there was a single linear order whose spectrum
contains all nonzero $\Delta^0_2$-degrees, but not $\bfz$.  It remains open
whether a linear order can have precisely the nonzero degrees in its spectrum.

For real closed subfields of $\Rz$, one might try to extend the above construction
to $\Delta^0_2$ sets $L$, but the method of $\Delta^0_2$-permitting does not
lend itself readily to the construction of computable Dedekind cuts.  In Theorem
\ref{thm:low}, we were able to compute the cuts in $F$, given finitely much information,
because the sequence of cuts at stages where $\sigma_m$ was active provides
a computable enumeration.  With $\Delta^0_2$-permitting, one would potentially
have to go back and forth between different approximations to the cuts, even
from one $\sigma_m$-stage to the next, and so it seems that the argument
for computability of the cuts would no longer hold.  We regard this, and also
the analogue of the result from \cite{M01}, as challenging questions.  Notice that,
while one can ask the same questions about spectra of real closed fields in general,
every archimedean real closed field with all nonzero c.e.\ degrees (or just a single minimal pair)
in its spectrum would have to be a subfield of $\Rz$.  So the question about subfields
of $\Rz$ is the essence of the more general question about all archimedean real closed fields.

\section{Nonarchimedean Real Closed Fields}
\label{sec:nonarch}

An ordered field $F$ is \emph{nonarchimedean} if it contains an element $x$
such that, for every $n\in\Z$, we have $n<x$ (with $\Z$ denoting the prime
subring of $F$).  Such an element $x$ is said to be \emph{positive infinite}.
More generally, $x$ is \emph{infinite} if either $x$ or $-x$ is positive infinite,
so the \emph{finite} elements are those $x$ such that, for some $n$, $-n<x<n$ in $F$.
An element is \emph{infinitesimal} if its reciprocal is infinite; $0$ itself is
usually also considered infinitesimal.
The finite elements form a local subring $R$ of $F$, with maximal ideal $\mathfrak{m}$
containing precisely the infinitesimal elements, and the archimedean
ordered field $R/\mathfrak{m}$ is called the \emph{residue field} $F_0$ of $F$.

Every element of $R$ realizes a Dedekind cut in the prime subfield $\Q$ of $F$.
Two distinct elements can realize the same cut, in which case they differ by an infinitesimal.
Each subfield of $\R$ is the residue field of many distinct ordered fields.
If the residue field is algebraic over $\Q$, then it is natural to think
of the algebraic closure of $\Q$ within $F$ as the (canonical) residue field
of $F$.  However, if $R$ contains an element $t$ transcendental over $\Q$,
then there is no canonical subfield of $R$ isomorphic to the residue field:
the Dedekind cut of $t$ will be realized by many elements of $R$, and there
is no natural way to choose just one of them to lie in a canonical subfield of $R$.
If $F$ is countable with domain $\omega$, then one might want to choose
the least element of $\omega$ which realizes this cut, and likewise for
each other cut.  However, there is no
reason to expect the set of those elements to form a subfield.  Moreover,
from the point of view of computability, it is only a $\Pi^0_1$ property for
two elements to realize the same Dedekind cut, and so, even if these least
elements realizing cuts did form a subfield, that subfield would not have a
natural presentation computable relative to the larger field $F$.
Our main theorem in this section, Theorem \ref{thm:nocompstd}, strengthens
this result by giving a computable nonarchimedean real closed field for which
the residue field has no computable presentation at all.  Before that, in Theorem
\ref{thm:prime}, we investigate one special case in which the residue field must have
a computable presentation.

The positive infinite elements of a nonarchimedean real closed field $F$ are partitioned into
the \emph{positive infinite multiplicative classes}, with two elements $x$ and $y$
lying in the same class if and only if there exists some $n\in\N$ for which
$y<x^n$ and $x<y^n$.  The order on the field gives a linear order $\L$ on these
positive infinite multiplicative classes, which we call the \emph{derived linear order}
of $F$.  One way to describe the field $F$ then is to enhance the language of fields
with constants $c_x$ for each element $x$ of the residue field and $d_A$ for
each positive infinite multiplicative class $A$, and to build a theory $T_F$
which includes the axioms of \RCF, the atomic diagram of the residue field
$(F_0,\set{ c_x}{x\in F_0})$ and the atomic
diagram of $\L$ on the elements $d_A$ (that is, for every $n>1$,
include either $n<d_A<d_A^n<d_B$ or $n<d_B<d_B^n<d_A$, depending
on whether $A<B$ in $\L$, so that $d_A$ and $d_B$ will be
positive infinite elements from distinct positive infinite classes).
Of course, the relation $x<y$ here abbreviates the statement that $(x-y)$ has no square root.

\label{prime}
The prime model of the theory $T_F$ is not necessarily $F$ itself, although
$F$ is a model of $T_F$.  To define the elements $c_x$ in $F$, noneffectively,
we go through the elements $x$ of $F_0$.  For each one in turn, the value
of $c_x$ in $F$ is determined if $x$ is algebraic over the previous ones;
if it is transcendental, then we can choose $c_x$ to be any element of $F$
from the residue class of $x$ in $F_0$.  The element $d_A$ in $F$ may be
any element from the class $A$.

The prime model of $T_F$ will be a real closed field whose residue
field is isomorphic to $F_0$ and whose derived linear order is isomorphic to $\L$,
and in which each of the positive infinite classes $A$ contains only elements built
from the $d_A$ and its predecessors by the field operations and the operations
for taking real closures.  Such a real closed field is therefore called a
\emph{prime-nonstandard field}.  In Section \ref{sec:questions}
we describe how a real closed field can fail to be prime-nonstandard.
For other questions about these fields, we refer the reader to \cite{O14},
which in turn uses unpublished work of Marker.

\begin{thm}
\label{thm:prime}
Fix any nonempty linear order $\L$ which has a left end point,
and any archimedean real closed field $F_0$.  Let $F$ be the
prime-nonstandard real closed field built from $\L$ and $F_0$. 
Then 
the spectrum of $F$ contains exactly those Turing degrees $\bfd\in\spec{F_0}$
such that $\bfd'\in\spec{\L}$.
\end{thm}
\begin{pf}
Both containments were proven in \cite{O14} for the case where $F_0$ is
the real closure of $\Q$:  in this case, every degree lies in $\spec{F_0}$
and can enumerate the cuts realized in $F_0$ and decide dependence there,
and so that proof showed that $\spec{F}=\set{\bfd}{\bfd'\in\spec{\L}}$.
This established the corollary that, for every linear order,
the preimage of its spectrum under the jump operation
is itself the spectrum of a real closed field.  Part of that
proof works equally well here for arbitrary $F_0$:  given a degree
$\bfd\in\spec{F_0}$ with $\bfd'\in\spec{\L}$, one builds
a $\bfd$-computable copy of $F$ in exactly the same manner,
using a $\bfd$-computable approximation of $\L$.
At stage $0$, for each $a\in\L$, we add one positive infinite element
$d_a$ to the field.  At certain stages the approximation may change its mind
and decides that, instead of having $a<b$, it wants $b<a$.
In this case, supposing that $a$ has higher priority than $b$,
we make $d_a$ and $d_b$ lie in the same positive infinite class
(by making $d_b=d_a^n$ for some large $n$)
add a new element to the field, redefine $d_b$ to be this new element,
and place the new $d_b$ to the left of $d_a$ in the field.
Each $d_a$ is redefined only finitely often by this finite-injury
process, so we have produced a copy of $F$.
This establishes the reverse containment.

The forward containment, that every $\bfd\in\spec{F}$
has $\bfd\in\spec{F_0}$ and $\bfd'\in\spec{\L}$, becomes apparent once one
realizes that when $\L$ has a left end point, the set of positive
infinite elements in a $\bfd$-computable copy of $F$
must be $\bfd$-decidable.  Indeed, fixing an element $y_0$
in the leftmost positive infinite multiplicative class
in a copy of $F$, we can enumerate the set of all positive
infinite elements in this copy:  they are those $x$ such that
$(\exists n\in\N)~y_0<x^n$.  However,
the positive finite elements are those $x>0$ for which
$(\exists n\in\N)~x < n$, and so the finiteness
of positive elements (hence of all elements) is $\bfd$-decidable.
This also gives decidability of the set of infinitesimals in $F$.
So, from any $\bfd$-computable copy of $F$,
both the subring of finite elements and the ideal
(within this subring) of infinitesimals are $\bfd$-decidable,
and their quotient is a $\bfd$-computable copy of $F_0$.

To see that $\bfd'\in\spec{\L}$, notice that $\L$ is interpretable
in $F$ by computable infinitary $\Sigma_2$ formulas.
The domain of the interpretation is the set of positive infinite elements
of $F$, modulo the relation of lying in the same multiplicative class,
and the linear order on this domain is lifted directly from $F$.
A $\bfd'$-oracle therefore allows us to recognize a positive infinite
element of $F$, to decide whether two such elements lie
in the same multiplicative class, and (if not) to decide
which one belongs to the larger class.  Since the multiplicative
classes of positive infinite elements form a copy of $\L$,
this yields a $\bfd'$-computable presentation of $\L$,
completing the proof.
\qed\end{pf}

In contrast to Theorem \ref{thm:prime}, the following is an example
of a computable nonarchimedean real closed field $F$ in which
the derived linear order $\L$ has no left end point
(and $F$ is the prime-nonstandard real closed field built
from $\L$ and the residue field $F_0$).  The theorem
shows that the conclusion of Theorem \ref{thm:prime}
need not hold in this case.  For its proof, it is important to
note that we continue here our convention of using right-leaning
Dedekind cuts.  (Later in this section, we will have to revise
this convention.)

\begin{thm}
\label{thm:nocompstd}
There exists a computable nonarchimedean real closed field $F$
whose residue field $F_0$ is not computably presentable.
Indeed, the spectrum of $F_0$ contains precisely the high Turing degrees.
\end{thm}
\begin{pf}
The construction of $F$ takes place on a tree $T$, using an
ordinary $\bfz''$ argument.  First, we line up our candidates
against which to diagonalize.  For each $e$, let $\P_{e}$
be the requirement that $W_e$ is infinite if and only
if the element $y_e$ (chosen from a sequence $y_0,y_1,\ldots$
of elements in our field $F$) realizes the same cut in $F$
as some rational number.  Understand that $y_e$ itself
will not be rational:  the question is whether there is a
rational from which $y_e$ differs by an infinitesimal.
For any degree $\bfd\in\spec{F_0}$, a $\bfd'$-oracle
will allow us to identify the element of (a $\bfd$-computable
copy of) $F_0$ realizing the same Dedekind cut as $y_e$,
and then to decide whether this element of $F_0$ is rational in $F_0$.
Thus the $\bfd'$-oracle will decide the $\Pi^0_2$-complete set $\Inf$,
proving the theorem.

In the tree $T$, all nodes at level $e$ will be $\P_e$-nodes,
devoted to satisfying $\P_e$, and the node on the true path at this level
will succeed in satisyfing it.  Each node $\alpha\in T$ will have two immediate
successors, labeled $\infty$ and $\fin$, with $\infty \prec\fin$,
representing the two outcomes $e\in\Inf$ and $e\in\Fin$.  As usual,
one node at each level $\leq s$ will be eligible at each stage $s+1$;
these nodes will all be comparable in $T$, and all nodes to their right in $T$
will be initialized at this stage.  The true path $P$ will contain the leftmost
node (under $\prec$) at each level which is eligible at infinitely many stages.

The strategy by which a $\P_e$-node $\alpha$ attempts to satisfy
$\P_e$ will involve choosing an element
$x_{\alpha,s}$ at each stage $s$.  This element will sit in between the 
infinitesimals and the standard (i.e., noninfinitesimal) positive elements,
without yet having been definitively assigned to either.  At some stage $s$
we may make $x_{\alpha,s}$ standard, by having some positive (standard)
integer $n$ in $F$ turn out to have $\frac1n < x_{\alpha,s}$; in this case we will choose
a new element $x_{\alpha,s+1} < x_{\alpha,s}$ which is in the new ``gap''
between infinitesimal and positive standard elements.  If there is some stage $s_0$
after which this never happens again (for this $e$ and $\alpha$), then
$x_{\alpha}=x_{\alpha,s_0}=\lim_s x_{\alpha,s}$ will be infinitesimal.
It will be a convenient feature of this construction
that, for all stages $s$ and all $\alpha,\beta\in T$,
we will have $x_{\alpha,s}<x_{\beta,s}$ in $F$ if and only if
$\alpha\prec\beta$ in $T$.  It follows that the same will
hold of $x_\alpha$ and $x_\beta$ whenever these limits exist.
(To be clear:  when $\alpha\subset\beta$, we define
$\beta\prec\alpha$ if $\alpha\widehat{~~}\infty\subseteq\beta$,
and $\alpha\prec\beta$ if $\alpha\widehat{~~}\fin\subseteq\beta$.)

The use of $x_{\alpha,s}$ is as follows.  The node $\alpha$
will have in mind a particular element $y_\alpha$ of $F$, which $\alpha$
describes during the construction by enumerating a cut in $\R$
to be filled by $y_\alpha$.  At each stage $s+1$, $\alpha$ will intend for
$y_\alpha$ to fill the cut of some rational number
$b_{\alpha,s}$.  The real intention, however, is that the cut of $y_\alpha$
should be the cut of a rational if $W_e$ turns out to be finite,
but not otherwise.  We will have $x_{\alpha,s}=b_{\alpha,s}-y_\alpha$,
so the uncertainty about whether $x_{\alpha,s}$ is infinitesimal
will mirror the uncertainty whether $y_\alpha$ lies in the cut of $b_{\alpha,s}$.
At stage $s+1$, if $W_{e,s}$ receives a new element,
then we make the cut of $y_\alpha$ in $F$ at this stage no longer
contain $b_{\alpha,s}$; instead, we redefine $b_{\alpha,s+1}\neq b_{\alpha,s}$
to be a distinct rational number, within the (new, narrower) cut of $y_\alpha$.
This ensures that $y_\alpha$ and $b_{\alpha,s}$ realize distinct cuts.

One possible outcome here is that $b_\alpha=\lim_s b_{\alpha,s}$ exists:
this occurs if $W_e$ is finite, and in this case $y_\alpha$ lies in the cut
of $b_\alpha$, so $\P_e$ is satisfied.  (This is the outcome $\fin$
of the node $\alpha$.)  Otherwise, we eliminate all rational numbers
from realizing the same cut as $y_\alpha$, thus showing that the cut
realized by $y_\alpha$ is irrational, so again $\P_e$ is satisfied.
(This is the outcome $\infty$ of the node $\alpha$.)  In either case,
$y_\alpha$ itself is definitely not rational; in fact it will
always be transcendental over $\Q$.  
Thus, the outcome $\infty$ corresponds to every $x_{\alpha,s}$
eventually being declared standard in $F$ (so that $b_{\alpha,s}$ and $y_\alpha$
do not realize the same cut).  The outcome $\fin$ corresponds to $y_\alpha$
staying in the same cut as the limit $b_\alpha$, in which case their difference
equals the limit $x_\alpha$, which exists and is infinitesimal.

Of course, when an element $x_{\alpha,s}$ is made standard at stage $s+1$, so is
every element $>x_{\alpha,s}$ then in $F$ (except those which are positive
and infinite).  So we will choose our elements 
$x_{\alpha,s}$ (for various $\alpha$) with this in mind, according to the guesses
by each node $\alpha$ about whether each node $\beta\subset\alpha$ will
act infinitely often or only finitely often.  Of course, $\alpha$ wants its own
$x_{\alpha,s}$ always to have the option of either becoming standard or staying infinitesimal.
If $\alpha$ thinks that $\beta$ will act infinitely often, then $\alpha$ picks
$x_{\alpha,s} < x_{\beta,s}$, so that $x_{\beta,s}$ can be made standard without
injuring $\alpha$'s strategy.  (Moreover, in this case, $x_{\alpha,s}$ can always be
made standard without injuring $\beta$'s strategy:  just wait until the next stage
at which $x_{\beta,s}$ is made standard, and make $x_{\alpha,s}$ standard
at the same stage.)  On the other hand, if $\alpha$ thinks that $\beta$
will act only finitely often -- and therefore, possibly, never again -- then it
chooses $x_{\alpha,s}>x_{\beta,s}$, so that $x_{\alpha,s}$ can be made standard
whenever needed without injuring the strategy of the higher-priority node $\beta$.

At stage $0$ we initialize all nodes, which means that we set all values
that vary over stages (that is, all values except the $y_\alpha$) to be undefined.
We make the first three elements of $F$ serve as the integers $0$, $1$, and $2$.
We fix the fourth element of $F$ to serve as $y_\lambda$, where $\lambda$ is the root node
of $T$, with $a_{\lambda,0}=1 < y_\lambda <b_{\lambda,0}=2$ in $F$.
(As usual, $a_{\lambda,s}$ and $b_{\lambda,s}$ denote the left and right
end points of the stage-$s$ approximation to the Dedekind cut of $y_{\lambda}$,
and likewise for other nodes.)
The set of all $y_\alpha$ (for all $\alpha\in T$ that are ever eligible)
will remain algebraically independent over $\Q$ at all stages.

At stage $s+1$, we first add one more element to $F$, according to
a systematic method of making $F$ a real closed ordered field.  The new element
may be a sum, product, or reciprocal of current elements, or a root
of a polynomial over current elements.  (What we have built at each
stage will always be just a finite fragment of $F$, necessarily not
closed under the field operations.)  For these purposes, the set of
elements $y_\alpha$ is treated as algebraically independent over $\Q$,
and each $b_{\alpha,s}$ currently defined is treated as infinitesimally
greater than the corresponding $y_\alpha$.

Notice that placing this new element within the current ordering of $F$
may require us to make a decision about where some $y_\alpha$ sits
inside its current interval in $F$.  This occurs if the new element is a
rational $q$ with $a_{\alpha,s}<q<b_{\alpha,s}$, of course, but it also
occurs if we add some high power $y_\alpha^j$ to $F$, for instance:
the placement of $y_\alpha^j$ in the rationals implicitly defines the placement
of $y_\alpha$ among the $j$-th roots of rationals.  In general, if the new
element is dependent on $y_\alpha$ (for one or more nodes $\alpha$),
then this may happen.  The assumption that each $y_\alpha$ is
infinitesimally less than $b_{\alpha,s}$ allows us to determine the
ordering on $F$.  
Moreover, we treat each element $x_{\alpha,s}$ as lying in its own
(positive infinitesimal) multiplicative class:  no new power of $x_{\alpha,s}$
will be smaller than any $x_{\beta,s}$ with $\beta\prec\alpha$.

The key here is that multiplicative-class nonequivalence
and infinitesimality (and also independence) are undefinable
by finitary formulas in the language of ordered fields, meaning
that no finite number of steps taken here can actually force any
difference $x_{\alpha,s}$ to be infinitesimal,
nor force any elements $x_{\alpha,s}$ and $x_{\beta,s}$ to lie in
distinct multiplicative classes.  (Nor can finitely many steps
make any $y_\alpha$ actually transcendental, although
in fact the elements $y_\alpha$ will all remain independent over $\Q$
at the end of the construction.)

Next, still at stage $s+1$, we proceed through the following steps for each eligible
node $\alpha$ at the level $e< s$ in $T$.  The root node $\lambda$ is eligible
at the beginning of the process.
\begin{itemize}
\item
If $y_\alpha$ is undefined, then we adjoin two new large consecutive
positive integers $a_{\alpha,s+1} < b_{\alpha,s+1}$ to $F$, along
with a new element $y_\alpha$ which is defined to lie between these
integers and to be transcendental over the fragment of $F$ defined so far.
(This value $y_\alpha$ will never change during the rest of the construction,
even if $\alpha$ is initialized.)
We make $\alpha\widehat{~~}\fin$ eligible at this stage.
\item
If $y_\alpha$ is defined but $\alpha$ has been initialized
since the last stage at which it was eligible, then we define
the elements $a_{\alpha,s+1}$ to be the greatest rational in
(the current finite fragment of) $F$ with $a_{\alpha,s+1}<y_\alpha$,
and adjoin a new rational $b_{\alpha,s+1}$ to $F$, between
$y_\alpha$ and the least element of $F$ greater than $a_{\alpha,s+1}$.
We also adjoin to $F$ the new element
$x_{\alpha,s+1}=b_{\alpha,s+1} - y_\alpha$
(since this difference is not yet defined), placing it so that, in the order on $F$,
we have $0 < x_{\alpha,s+1} <\frac1n$ for all positive integers $n$
already in $F$, and 
$$x_{\alpha,s+1} < x_{\beta,s+1} \iff\alpha\prec\beta$$
for all $\beta$ with $x_{\beta,s}$ defined.
Notice that, although this $b_{\alpha,s+1}$ is
the right end point of the current cut of $y_\alpha$ in $F$.
this does not preclude $y_\alpha$ from realizing the same
Dedekind cut in $F$ as $b_{\alpha,s+1}$, since the difference
$x_{\alpha,s+1}$ between them could be infinitesimal.
We make $\alpha\widehat{~~}\fin$ eligible at this stage.
\item
Otherwise, let $s' < s$ be the most recent stage at which $\alpha$
was eligible.
If $W_{e,s}=W_{e,s'}$, then we deem $\alpha$ \emph{inactive}
at stage $s+1$ and make $\alpha\widehat{~~}\fin$ eligible.
If $W_{e,s}\neq W_{e,s'}$, then we deem $\alpha$ to be
\emph{active} at this stage, meaning that we will include $\alpha$
in the action taken as we complete stage $s+1$, after all substages
are finished.  In this case, $\alpha\widehat{~~}\infty$ becomes eligible.
\end{itemize}

After completing this step for the eligible node at level $s-1$,
we let $A_s$ be the set of those $\alpha$ which are active at this stage.
Recall that by definition $x_{\alpha,s} = b_{\alpha,s}-y_\alpha$ in $F$.
By induction, we know, for all $\alpha,\beta\in T$, that
$x_{\alpha,s} < x_{\beta,s}$ in $F$ if and only if $\alpha\prec\beta$.
We add to $F$ a new (finite) large positive integer $n$,
placed in the order so that, for all $\beta\in T$,
$$\frac1n < x_{\beta,s} \iff (\exists \alpha\preceq\beta)\alpha\in A_s.$$
The ordering $\prec$ on $T$, along with the inductive assumption above,
makes this consistent with \RCF.  Now that
we have $x_{\alpha,s}>\frac1n$ for these $\alpha$,
we immediately (for each of these $\alpha$) choose one new rational number
$b_{\alpha,s+1}$ from the interval $(a_\alpha,b_{\alpha,s})$ in $F$.
In particular, we go through the $A_s$ in order under $\subseteq$, starting
with the shortest $\alpha$ there.  For each $\alpha$ in its turn, having
determined the appropriate interval, we choose the new $b_{\alpha,s+1}$
to be the first rational we find in that interval, under some fixed ordering of $\Q$.
We then choose each $a_{\alpha,s+1}$ between $a_{\alpha,s}$ and
$b_{\alpha,s+1}$ in such a way that none of the first $s$ cuts
in an enumeration of the Dedekind cuts in the real closure of
$\Q(y_{\beta,s+1}~:~\beta\widehat{~~}\infty\subseteq\alpha)$
lies in the open interval $(a_{\alpha,s+1},b_{\alpha,s+1})$.
This interval is the next step in our computation of the Dedekind cut of $y_\alpha$:
we now define $y_\alpha$ to lie in this interval, and think of it as
(without formally defining it to be) infinitesimally close to $b_{\alpha,s+1}$.
The new element $x_{\alpha,s+1}=b_{\alpha,s+1} - y_\alpha$
is also defined right now for each $\alpha\in A_s$:  it is less than those positive
elements of $F$ that are already standard, but greater than all positive elements of $F$
which have not yet been made standard, and these new elements
$x_{\alpha,s+1}$ lie in distinct positive-infinitesimal multiplicative classes,
ordered so as to obey the rule that $x_{\alpha,s+1} < x_{\beta,s+1}$ in $F$
iff $\alpha\prec\beta$ in $T$.

Of course, for every $\alpha\in A_s$ and every $\beta\succ\alpha$ on $T$,
this process also makes $x_{\beta,s}>\frac1n$.  All such $\beta$ are initialized
at this stage, and are said to have been \emph{injured} by this process.
(Notice that, if $A_s=\emptyset$, then each eligible node at this stage
is the rightmost node at its level, and so no initialization takes place
anywhere.)  On the other hand, those $\beta\in T$ such that
$(\forall\alpha\in A_s)\beta\not\succeq\alpha$ must have $x_{\beta,s}<\frac1n$,
by our placement of $n$ (or else $x_{\beta,s}$ and $b_{\beta,s}$ are undefined),
and so those $x_{\beta,s}$ are not forced into the standard part of $F$ by this step,
although they could yet become standard at a future stage.  For those $\beta$,
we keep $x_{\beta,s+1}=x_{\beta,s}$ and $b_{\beta,s+1}=b_{\beta,s}$,
but choose $a_{\alpha,s+1}$ to be the greatest rational in $F_s$
lying in the interval
$[a_{\alpha,s},b_{\alpha,s})$.  (This ensures that we do build Dedekind cuts
even for those $y_\beta$ with $\beta$ to the left of the true path.)

This completes the stage $s+1$,
and we define $F$ to be structure built over all stages by this construction.
Since every step is consistent with the theory \RCF, the initial steps
taken at all stages collectively ensure that $F$ is a
computable real closed field, although not archimedean.
It is clear that the true path $P$ through the tree $T$ exists
(since at each stage $s$, some node at each level $< s$ was eligible),
and that in fact the path $P$ is exactly the set $\Fin$.

For each node $\alpha\subset P$, say of length $e$,
there is a stage $s_0$ after which $\alpha$ is never again
initialized.  From then on, each time $\alpha$ is eligible,
we check whether $W_e$ has changed since the last time.
If this only happens finitely often, then $b_\alpha=\lim_s b_{\alpha,s}$
exists, and $y_\alpha$ lies in the same Dedekind cut as $b_\alpha$,
since for every rational $q<b_\alpha$ subsequently entering $F$,
we made $y_\alpha>q$.  So in this case $y_e=y_\alpha$ satisfies $\P_e$.
On the other hand, if $W_e$ is infinite, then $b_{\alpha,s}$
was redefined at infinitely many stages $s$.  Each time,
the new $b_{\alpha,s}$ was the first available rational
satisfying the order requirements.  It follows that
$y_\alpha$ cannot realize the same cut as any rational in $F$,
and so again $\P_e$ is satisfied by $y_e=y_\alpha$.
Thus the sequence $y_0,y_1,\ldots$ defined here
is the sequence we described at the start.  (It is not
a computable sequence, since it is defined using the true path $P$,
but it will serve our purposes below.)

Now suppose that a degree $\bfd$ can compute a copy $E_0$
of $F_0$, the residue field of $F$.  Then, with a $\bfd'$ oracle,
we could identify the unique element $z_0$ of $E_0$
which realizes the same Dedekind cut that $y_0=y_\lambda$
realizes in $F$ (where $\lambda$ is the root of $T$).  Moreover,
since $\bfd$ can enumerate the rationals in $E_0$,
the $\bfd'$ oracle lets us decide whether $z_0$ is rational in $E_0$,
and hence whether $y_0$ realizes a cut in $F$ whose right half
has a least element.  (That is, we have decided whether the cut of $y_0$
is the cut of some rational.)  The element $y_0$ of $F$ is irrational
there in any case, indeed transcendental, but its cut realizes a rational
number if and only if $W_0$ was finite, so we have decided finiteness of $W_0$.

Having done so, we know which of the two level-$1$ nodes in $T$
lies on the true path.  Let $\alpha_1$ be this node, and set $y_1=y_{\alpha_1}\in F$.
Then by the same process, we can use our $\bfd'$-oracle to decide whether
or not $W_1$ is finite, and we can continue recursively and decide
the $\Pi^0_2$-complete set $\Inf$ using just the $\bfd'$ oracle.
Thus, every degree $\bfd\in\spec{F_0}$ must be high.

Conversely, if $\bfd$ is high, then $\bfd$ satisfies the demands
of Theorem \ref{thm:spec} for lying in the spectrum of $F_0$.
In particular, $F_0$ is algebraic over the set of all elements
$y_\alpha=\lim_s y_{\alpha,s}$ (since whenever a node $\alpha$
was initialized, its old $y_{\alpha,s}$ became rational), and hence
we have a transcendence basis consisting of those $y_\alpha$
which do not lie in the cut of an algebraic number.  (To be clear:
$F_0$ consists of classes of standard elements of $F$, modulo
infinitesimals, and the collection of classes of such $y_\alpha$
forms a transcendence basis for $F_0$.)  So we wish to enumerate
the cuts of those $y_\alpha$ for which $\alpha\hat{~}\infty$
is on the true path $P$; all other nodes $\alpha$ either were initialized
infinitely often, or had $\alpha\hat{~}\fin\subset P$,
or were eligible only finitely often, and in the latter two cases
$y_\alpha$ wound up in the cut of an algebraic number.
But since $\bfd'\geq\bfz''$, a $\bfd$-oracle can approximate
the true path.  Let $\alpha_0,\alpha_1,\ldots$ list all
nodes on the tree $T$.  When $\bfd$ first thinks that
$\alpha_0\hat{~}\infty\subset P$ and that $\alpha_0$ is never initialized
after stage $s$, it begins enumerating a cut $C_0$ for $y_{\alpha_0,s}$.
If the $\bfd$-approximation later changes its mind, it stops
using this cut for $\alpha_0$ and instead uses it for the
$y_{\alpha_n,t}$ currently being enumerated by $C_1$,
adjusting $y_{\alpha_n,t}$ by adding
a rational as needed to make it lie in the cut $C_0$ as
enumerated up till now.  ($C_1$ in turn takes over the $y_{\alpha_m,t}$
currently being enumerated by $C_2$, plus a rational, and so on.)
If, still later on, the $\bfd$-approximation
decides that $\alpha_0\hat{~}\infty$ is on $P$ after all,
it switches $C_0$ to enumerate the cut of the new $y_{\alpha_0,s}$,
plus a rational, and moves the current $y_{\alpha_n,t}$ back to $C_1$,
with all later cuts shifting similarly.  Since the $\bfd$-approximation
actually does converge to the characteristic function of $P$,
we wind up with an enumeration of cuts in which each cut
contains some $(y_\alpha+q_\alpha)$ for some $\alpha\hat{~}\infty\subset P$
and some $q_\alpha\in\Q$, and every $\alpha\hat{~}\infty\subset P$
has its $y_\alpha$ in exactly one of these cuts, up to a rational difference.
Thus we have satisfied Condition (3) of Theorem \ref{thm:spec},
so $\bfd\in\spec{F_0}$, and $\spec{F_0}$ contains precisely the high degrees.
\qed\end{pf}

\begin{cor}
\label{cor:notenough}
The spectrum of a prime-nonstandard real closed field $F$ is not determined by the spectra
of its residue field $F_0$ and its derived linear order $\L$.
\end{cor}
\begin{pf}
Theorem \ref{thm:nocompstd} yields a computable $F$ for which
$\spec{F_0}$ contains exactly the high degrees.  Let $\L$ be the derived
linear order for this $F$ (which in fact is the order $\omega^*$ of the
negative integers), and use the construction from \cite{O14} to build
another real closed field $E$ with the same derived order, but with
the field $\Rz$ of computable real numbers as its residue field.
Now $\spec{\Rz}=\spec{F_0}$, and the derived orders are isomorphic,
but we claim that $E$ has no computable presentation (and hence
that $\spec{E}\neq\spec{F}$).

Being nonarchimedean, $E$ does not quite allow us to enumerate
its right-leaning cuts, for the same reason that $F$ did not (in Theorem
\ref{thm:nocompstd}).  However, we can use a presentation of $E$
to enumerate all non-strict Dedekind cuts of real numbers
in $[0,1]$ realized in the residue field $\Rz$:  for each $y\in E$
with $0\leq y\leq 1$ in $E$,
build both the left-leaning cut $((-\infty,y], (y,+\infty))$ of $y$ 
and the right-leaning cut $((-\infty,y), [y,+\infty))$ of $y$ in $\Q$.
Notice that if $y=q+\ep$ for a rational $q$ and positive infinitesimal $\ep$,
then both of these equal the left-leaning cut of $q$; whereas if $y=q-\ep$,
then both equal $q$'s right-leaning cut.  So we have satisfied precisely
the first item in Proposition \ref{prop:bridge}, and therefore
the presentation of $E$ must be of high degree.
\qed\end{pf}

We salvage the following result about residue fields, which, in light of
Theorem \ref{thm:nocompstd}, is the best possible statement of its form.
\begin{prop}
\label{prop:residue}
If $F$ is a real closed field of Turing degree $\bfc$, then its residue
field $F_0$ has a presentation in every degree $\bfd\geq\bfc$ with $\bfd'\geq\bfc''$.
\end{prop}
\begin{pf}
Let $\bfd$ be such a degree. Since $\bfc$ can enumerate the Dedekind cuts
realized in $F_0$ (by waiting until an element of $x\in F$ satisfies
$(\exists n\in\N)~-n<x<n$, and then enumerating the cut of $x$),
so can $\bfd$.  Moreover, since $\bfc''$ can decide whether the real
number realized by one of these cuts is transcendental over the reals
realized by the preceding cuts, it can decide the dependence relation
on these cuts, and therefore so can $\bfd'$.  Thus $\bfd$ satisfies
Condition (2) of Theorem \ref{thm:spec} for lying in $\spec{F_0}$.
\qed\end{pf}

Theorem \ref{thm:nocompstd} recalls the related result
\cite[Theorem 4.2]{KL13}.  There, Knight and Lange built
a computable real closed field $F$ for which no residue field section
can be a $\Sigma^0_2$ subset of $F$.  (A \emph{residue field section}
of $F$ is an archimedean real closed subfield of $F$
containing exactly one element from each Dedekind cut realized
by a finite element of $F$.  They had previously shown that such an $F$
must have a $\Pi^0_2$ residue field section.)  We see no direct proof
of either of Theorem \ref{thm:nocompstd} or \cite[Theorem 4.2]{KL13}
from the other.  However, although the two constructions were developed
independently, they bear a close resemblance, and any further work
extending either of them might well yield results extending the other as well.

\section{Further Results and Questions}
\label{sec:questions}

Real closed fields in general, and even archimedean ones, do realize many
more spectra than linear orders:  every upper cone, say above $\bfc$,
is the spectrum of the real closure of $\Q(x)$, where $x$ is a real number
of Turing degree $\bfc$, and every $\bfc'$ is the jump degree of $\Rc$.
We do not yet know of any spectrum of a structure
which cannot be realized as the spectrum of a real closed field.  However,
we conjecture that spectra such as the set of all non-low$_n$ degrees may
be beyond the reach of real closed fields, especially archimedean ones.

It should be noted that the method used so far to show completeness
for spectra of various classes of structures (as described in Section \ref{sec:intro})
will not suffice for real closed fields.  This method can be described
as the construction of a \emph{computable functor} with a computable inverse,
usually between an arbitrary graph and a member of the desired class;
see \cite{MPSS16} for details and definitions,
and \cite{HTMMM16} for an alternative version using effective interpretations
of one structure in another.  However, such transformations preserve many other properties.
In particular, they preserve the automorphism group:  if we have computable
functors (as in \cite[Defn.\ 1.8]{HTMMM16}) between $\A$ and $\B$,
then the automorphism groups of these structures must be isomorphic.
However, since a real closed field $F$ always has an underlying linear order $<$
(definable in the field, whether or not $<$ is included in the signature),
a non-trivial automorphism $h$ of $F$ must be a non-trivial automorphism
of this order, and therefore can never have $h^n(x)=x$:  if $x<h(x)$, then
$x<h(x)<h(h(x))<\cdots$, and similarly if $h(x)<x$.  On the other hand,
many graphs have non-identity automorphisms of finite order,
and so each such graph cannot be computably bi-transformable with any
real closed field.  Thus, if a method is to be found for showing that every
spectrum of a graph is the spectrum of a real closed field, it will have
to be a new method, preserving spectra while failing to preserve other standard
computable-model-theoretic properties.

The situation is made more murky by the existence of non-archimedean
real closed fields which fail to be prime, in the sense described in Section
\ref{sec:nonarch} (page \pageref{prime}).
For a simple example, start with the real closure of $\Q$ as the residue field
and the derived linear order $\L$ with just one element:  we get a prime-nonstandard
real closed field $F$ from these.  Fix a positive infinite element $x\in F$,
and now adjoin an element $y=x^{\sqrt2}$ to $F$.  That is, $y$ should be
transcendental over $F$ and satisfy
$$ x^q<y \iff q<\sqrt2,$$
with the obvious extension of this $<$ to all other elements of $F$.
This $F$ is \emph{non-prime}, in the sense that it is not the prime
model of the theory which includes $\RCF$ along with the existence of one positive
infinite element (i.e, the theory saying that this $\L$ embeds into
the derived linear order).  In fact, this $F$ still has a computable presentation,
but more complicated non-prime nonstandard real closed fields,
with elements of the form $x^a$ for noncomputable real numbers $a$,
may realize new spectra. 
Section \ref{sec:nonarch} represents some
early steps towards addressing the question of spectra of nonarchimedean
real closed fields, but it seems clear that much more remains to be done.

\parbox{4.7in}{
{\sc
\noindent
Department of Mathematics \hfill \\
\hspace*{.1in}  Queens College -- C.U.N.Y. \hfill \\
\hspace*{.2in}  65-30 Kissena Blvd. \hfill \\
\hspace*{.3in}  Flushing, New York  11367 U.S.A. \hfill \\
Ph.D. Programs in Mathematics \& Computer Science \hfill \\
\hspace*{.1in}  C.U.N.Y.\ Graduate Center\hfill \\
\hspace*{.2in}  365 Fifth Avenue \hfill \\
\hspace*{.3in}  New York, New York  10016 U.S.A. \hfill}\\
\medskip
\hspace*{.045in} {\it E-mail: }
\texttt{Russell.Miller\at {qc.cuny.edu} }\hfill \\
}\\

\parbox{4.7in}{
{\sc
\noindent
University of Puerto Rico at Mayaguez \hfill \\
\hspace*{.1in}  Department of Mathematical Sciences \hfill \\
\hspace*{.2in}  Call Box 9000 \hfill \\
\hspace*{.3in}  Mayaguez, PR 00681-9018  U.S.A. \hfill}\\
\medskip
\hspace*{.045in} {\it E-mail: }
\texttt{victor.ocasio1\at {upr.edu} }\hfill \\
}

\end{document}